\documentclass[12pt,letterpaper]{article}

\usepackage{amssymb,amsfonts,amscd,amsthm}
\usepackage{thmtools}
\usepackage[all,arc]{xy}
\usepackage{enumerate}
\usepackage{mathrsfs}
\usepackage{tikz}
\usepackage[left=0.9in,top=0.9in,right=0.9in,bottom=0.9in]{geometry}
\usepackage{mathtools}
\usepackage{hyperref}
\usepackage{cleveref}
\usepackage{color}%colors
\usepackage{graphicx}
\usepackage{enumitem} %Allows resuming enumeration and numbering by section
\graphicspath{ {TexImages/} }

\newtheorem{thm}{Theorem}[section]
\newtheorem{cor}[thm]{Corollary}

\newtheorem{lem}[thm]{Lemma}
\newtheorem{claim}[thm]{Claim}

\newtheorem{conj}[thm]{Conjecture}

\theoremstyle{definition}
\newtheorem{defn}{Definition}

\hypersetup{
	colorlinks,
	citecolor=blue,
	filecolor=blue,
	linkcolor=blue,
	urlcolor=blue,
	linktocpage
}

\setlist[enumerate]{itemsep=2ex, topsep=2ex} %spaces out enumerate/itemize better
\setlist[itemize]{itemsep=2ex, topsep=2ex}

\newcommand{\al}{\alpha}
\newcommand{\be}{\beta}
\newcommand{\gam}{\gamma}

\newcommand{\ep}{\varepsilon}

\newcommand{\Om}{\Omega}

\newcommand{\half}{\frac{1}{2}}

\newcommand{\sm}{\setminus}
\newcommand{\sub}{\subseteq}

\renewcommand{\SS}[1]{\textcolor{red}{#1}}

\newcommand{\dist}{\mathrm{dist}}

\title{Equivalence of Small and Large Quasi-kernels}
\author{Sam Spiro\footnote{Dept.\ of Mathematics and Statistics, Georgia State University, sspiro@gsu.edu}.}
\date{\today}

\begin{document}
	
	\maketitle
\begin{abstract}
	A quasi-kernel of a digraph $D$ is an independent set which can reach every vertex of $D$ by a directed path of length at most 2.  The Small Quasi-kernel Conjecture states that every $n$-vertex source-free digraph $D$ contains a quasi-kernel $Q$ with $|Q|\le n/2$, and the Large Quasi-kernel Conjecture states that every $n$-vertex digraph $D$ contains a quasi-kernel $Q$ with $|N^+[Q]|\ge n/2$.  Our main result shows for any $\frac{1}{2}<\beta<1$ an equivalence between the following two statements:
	\[\textrm{Every }n\textrm{-vertex source-free digraph contains a quasi-kernel with } |Q|\le n-\Omega(n^\beta).\]
	\[\textrm{Every }n\textrm{-vertex digraph contains a quasi-kernel with } |N^+[Q]|\ge \Omega(n^{2-\frac{1}{\be}}).\]
	That is, any sublinear bound for one problem is equivalent to a corresponding sublinear bound for the other.  This combined with previous work establishing an equivalence for linear bounds shows that the two problems are largely equivalent to each other.
	
	We additionally give a slight refinement to the recent breakthrough of Rutschmann for these problems which, together with the result above, gives the best known bounds of $n-\Om(n^{3/4})$ and $\Omega(n^{2/3})$ for the Small and Large Quasi-kernel Conjectures, respectively.  Our proofs are due to ChatGPT's Astra 6; the writing is all our own.  %This article is not intended for publication.
\end{abstract}

\section{Introduction}
This paper concerns problems about (finite) digraphs, for which we recall a few definitions.  Given vertices $u,v$ of a digraph $D$, we define $\dist(u,v)$ to be the length of a shortest directed path from $u$ to $v$, and for a set of vertices $S$ we define $\dist(S,v)=\min_{u\in S} \dist(u,v)$.   We say a set of vertices $S$ is \textit{independent} if there exist no arcs between two vertices of $S$.  The main definition for this paper is the following.

\begin{defn}
	A set of vertices $Q$ of a digraph $D$ is said to be a \textit{quasi-kernel} if (1) $Q$ is an independent set and (2) $\dist(Q,x)\le 2$ for all $x\in V(D)$ (that is, every vertex $x\in V(D)$ can be reached from $Q$ by a directed path of length at most 2).
\end{defn}
It was shown by Chv\'atal and Lov\'asz \cite{chvatal2006every} that every digraph has a quasi-kernel.  Because of this, it is natural to ask extremal problems about the set of quasi-kernels for a given digraph.  One conjecture in this direction is the Small Quasi-kernel Conjecture due to P.L.\ Erd\H{o}s and Sz\'ekely from 1976 (though it only formally appeared later on \cite{erdos2010two}).  Recall that a digraph is \textit{source-free} if it contains no vertices with in-degree 0. 

\begin{conj}[Small Quasi-kernel Conjecture]\label{conj:smallquasi}
	Every $n$-vertex source-free digraph $D$ contains a quasi-kernel $Q$ with $|Q|\le\half n$.
\end{conj}
We recently posed a variant of this problem which asks for the existence of ``large'' quasi-kernels.  For this we recall that for a digraph $D$ and a set $S\sub V(D)$ we write
\[N^+_D(S)=\{v\in V(D)\sm S:\exists u\in S,\ uv\in E(D)\},\qquad N^+_D[S]=S\cup N^+_D(S),\]
and we denote these sets simply by $N^+(S)$ and $N^+[S]$ whenever $D$ is clear from context.
\begin{conj}[Large Quasi-kernel Conjecture]
	Every $n$-vertex digraph $D$ contains a quasi-kernel $Q$ with $|N^+[Q]|\ge \half n$.
\end{conj}
The Small and Large Quasi-kernel Conjectures are equivalent in a strong sense: for all $0<\ep\le \half$, proving every $n$-vertex source-free digraph has a quasi-kernel $Q$ with $|Q|\le (1-\ep) n$ is equivalent to proving every $n$-vertex digraph has a quasi-kernel $Q$ with $|N^+[Q]|\ge \ep n$.  The forward direction of this was shown in \cite{spiro2026generalized}, the reverse was shown recently by Ai, Gutin, Liu, and Meng~\cite{ai2026small} with a somewhat weaker relationship in the constants being proven earlier by Ai, Liu, and Peng~\cite{ai2024variable}.

The Small Quasi-kernel Conjecture has been solved for various classes of digraphs, and we refer the reader to the nice survey by P.L.\ Erd\H{o}s, Gy\H{o}ri, Mezei, Salia, and Tyomkyn \cite{erdHos2025small} for a more thorough overview of these and related results.

Our focus in this paper is in proving bounds for quasi-kernels which hold for arbitrary digraphs.  In recent breakthrough work, Rutschmann~\cite{rutschmann2026small} improved the best known lower bound for the Large Quasi-kernel Conjecture to $n^{2/3-o(1)}$ for all $n$-vertex digraphs.  Additionally, he proved a general transference theorem from the Large Quasi-kernel Conjecture to the Small Quasi-kernel Conjecture, thereby establishing a corresponding upper bound of the form $n-n^{2/3-o(1)}$ for the Small Quasi-kernel Conjecture as well.

The main result of this paper is a stronger transference theorem between the conjectures which shows that establishing sublinear bounds for one problem is equivalent to establishing a corresponding sublinear bound for the other.

\begin{thm}\label{equivalence}
	The following statements are equivalent for any real number $0<\al<1$:
	\begin{itemize}
		\item[(a)] There exists some $c>0$ such that every $n$-vertex digraph has a quasi-kernel $Q$ with \[|N^+[Q]|\ge c n^\al.\]
		\item[(b)] There exists some $c'>0$ such that every $n$-vertex source-free digraph has a quasi-kernel $Q$ with \[|Q|\le n-c'n^\frac{1}{2-\al}.\]
	\end{itemize}
\end{thm}
Substituting $\be=\frac{1}{2-\al}$ recovers the version of this result stated in the abstract.

\Cref{equivalence} combined with the Large Quasi-kernel bound of Rutschmann gives an improved Small Quasi-kernel bound of $n-n^{3/4-o(1)}$.  The second result of our paper slightly refines Rutschmann's argument to get rid of the $o(1)$ terms in both exponents.

\begin{thm}\label{refinemnet}
	Every $n$-vertex digraph has a quasi-kernel $Q$ with
	\[|N^+[Q]|\ge \frac{1}{64} n^{2/3},\]
	and every $n$-vertex source-free digraph has a quasi-kernel $Q$ with
	\[|Q|\le n-\frac{1}{768}n^{3/4}.\]
\end{thm}

The constants of this theorem are far from optimized.

\iffalse
This together with a more quantitative version of \Cref{equivalence} gives a corresponding bound for the Small Quasi-kernel Conjecture.

\begin{cor}
	\SS{Insert}
\end{cor}
We note that the constant here and in \Cref{refinemnet} are not optimized.\fi 

\textbf{AI Declaration}.  Every proof in this paper is due to ChatGPT's Astra 6, though most of its ideas appeared previously in entirely human generated papers \cite{ai2024variable,rutschmann2026small,spiro2026generalized}.  Notably, we prompted Astra several times to go beyond the bound of $\Om(n^{1/2})$ for the Large Quasi-kernel Conjecture, but it failed to do so until the paper of Rutschmann's appeared.  While \Cref{refinemnet} gives the best known bounds for these problems, we emphasize that the main ideas all originate with Rutschmann.   Every word of this paper was written by the author and we take full responsibility for its mathematical contents.

\textbf{Acknowledgments}.  Access to ChatGPT Astra was provided to us by the ChatGPT for Academic Researchers program.  We thank Peter L.\ Erd\H{o}s for helpful comments.

\section{Preliminaries}
In what follows, we recall the basic digraph notation $D-S$ where $S\sub V(D)$ to mean the digraph obtained by deleting every vertex of $S$ from $D$, as well as $D[S]$ to mean the digraph obtained by deleting every vertex not in $S$ from $D$.

Our proof of \Cref{refinemnet} as well as part of our proof of \Cref{equivalence} requires us to construct quasi-kernels satisfying certain conditions. This is accomplished through the following lemma, which is a straightforward adaptation of the original proof by Chv\'atal and Lov\'asz~\cite{chvatal2006every} showing quasi-kernels exist.
\begin{lem}\label{CL-move}
	Let $D$ be a digraph and $I$ an independent set of $D$.  For every quasi-kernel $Q'$ of $D':=D-N^+_D[I]$, there exists a quasi-kernel $Q$ satisfying 
	\[Q'\sub Q,\]
	and 
	\[I\sub N^+_D[Q].\]
\end{lem}
\begin{proof}
	Take $Q=Q'\cup (I\sm N_D^+(Q'))$.  Since $Q'\sub V(D)-N^+_D[I]$ the set $Q$ is independent by construction.  The set $Q'\sub Q$ reaches every vertex of $V(D')=V(D)-N^+_D[I]$ within two steps, and $Q$ by definition reaches every vertex of $I$ within one step and hence every vertex of $N^+_D[I]$ within two steps.  This proves $Q$ is a quasi-kernel, and by construction it satisfies all of the conditions of the lemma.
\end{proof}
Iteratively applying this lemma gives the following.
\begin{cor}\label{CL-sequence}
	Let $D_1,\ldots,D_t$ be a sequence of digraphs obtained by having $D_{i+1}=D_i-N^+_{D_i}[I_i]$ for some independent set $I_i$ of $D_i$.
	\begin{itemize}
		\item[(a)] For every quasi-kernel $Q_i$ of $D_i$, there exists a quasi-kernel $Q_1$ of $D_1$ with $Q_i\sub Q_1$.
		\item[(b)] There exists a quasi-kernel $Q_1$ of $D_1$ with $I_1\cup \cdots \cup I_{t-1} \sub N^+_{D_1}[Q_1]$.
	\end{itemize}
\end{cor}
\begin{proof}
	We prove (a) by induction on $i$.  The $i=1$ case is immediate.  For $i>1$ assume we have proven the result for $i-1$.  By \Cref{CL-move} there exists a quasi-kernel $Q_{i-1}$ of $D_{i-1}$ containing $Q_i$, and applying the inductive hypothesis to $Q_{i-1}$ gives a quasi-kernel of $D_1$ containing $Q_{i-1}\supseteq Q_i$.
	
	We prove (b) via showing by reverse induction on $1\le i\le t$ that there exists a quasi-kernel $Q_i$ of $D_i$ with $I_i\cup I_{i+1}\cdots \cup I_{t-1}\sub N^+_{D_i}[Q_i]$.  The $i=t$ case follows by taking any quasi-kernel of $D_t$ (which exist either by Chv\'atal and Lov\'asz~\cite{chvatal2006every} or by applying \Cref{CL-move}). For $i<t$ assume we have already inductively constructed $Q_{i+1}$.  We can then apply \Cref{CL-move} with respect to $Q'=Q_{i+1}$ and $I=I_i$ to obtain a quasi-kernel $Q_i$ of $D_i$ which satisfies
	\[N^+_{D_i}[Q_i]\supseteq I_i\cup N_{D_{i+1}}^+[Q_{i+1}]\supseteq I_i\cup I_{i+1}\cdots \cup I_{t-1}.\]
	This establishes the inductive step, proving the result.
\end{proof}
We note that \Cref{CL-sequence} more or less appears in the work of Rutschmann~\cite{rutschmann2026small} with the digraphs $D_i$ being referred to as CL-subgraphs due to its inspiration from the proof of Chv\'atal and Lov\'asz.

In addition to these facts about digraphs, we will frequently use without comment the fact that the function $f(x)=x^\gam$ for $0< \gam\le 1$ is concave, and hence
\[x^\gam-y^\gam\ge \gam x^{\gam-1}(x-y)\]
for all $x>0$ and $y\ge 0$.

\section{Proof of \Cref{equivalence}}
We begin by proving small quasi-kernels imply large quasi-kernels.  The argument is essentially a more careful analysis of our proof from \cite{spiro2026generalized} showing the Small Quasi-kernel Conjecture implies the Large Quasi-kernel Conjecture.

\begin{lem}\label{small implies large}
	If $0<c<1$ and $\half<\be<1$ are real numbers such that every $n$-vertex source-free digraph $D'$ has a quasi-kernel $Q'$ with $|Q'|\le n-c n^\be$; then every $n$-vertex digraph $D$ has a quasi-kernel $Q$ with $|N^+[Q]|\ge \frac{1}{8} c^{1/\be} n^{2-\frac{1}{\be}}$.
\end{lem}
\begin{proof}
	Assume for contradiction that the result is false and choose a counterexample $D$ with $n$ as small as possible.
	\begin{claim}
		The digraph $D$ does not contain sources.
	\end{claim}
	\begin{proof}
		Assume $D$ had a source $v$ and let $d:=|N^+_D[v]|$.  Because $D$ is a minimal counterexample, there exists a quasi-kernel $Q'$ of $D':=D-N^+_D[v]$ with $|N^+_{D'}[Q']|\ge  \frac{1}{8} c^{1/\be} (n-d)^{2-1/\be}$.  Observe that $Q:=Q'\cup \{v\}$ is a quasi-kernel of $D$ since $v$ is a source and $N^+_D[v]\cap Q'=\emptyset$.  In total,
		\begin{align*}|N^+_D[Q]|\ge |N^+_{D'}[Q']|+|N^+_D[v]|&\ge \frac{1}{8} c^{1/\be}(n-d)^{2-1/\be}+d\\ &\ge \frac{1}{8} c^{1/\be}[(n-d)^{2-1/\be}+d] \ge \frac{1}{8} c^{1/\be} n^{2-1/\be},\end{align*}
		with this last inequality using that $2-1/\be\le 1$.  This contradicts $D$ being a counterexample, proving the result.
	\end{proof}
	Let $m$ be the maximum size of $|N^+_D[Q]|$ amongst all quasi-kernels $Q$ of $D$, and with some foresight we take
	\[k:=\left\lfloor \frac{n}{m}\right\rfloor\ge 1.\]
	Now define a digraph $D'$ by taking each vertex $v\in V(D)$ and adding $k$ new vertices $v_1,\ldots,v_k$ together with the arcs $vv_i$ for all $i$; that is, we form $D'$ by attaching $k$ directed leaves to each vertex in the original digraph $D$.  Let $Q'$ be a smallest quasikernel of $D'$, noting by hypothesis that 
	\[|Q'|\le |V(D')|-c|V(D')|^\be=(k+1)n-c((k+1)n)^\be.\]
	
	It is not difficult to verify that $Q:=Q'\cap V(D)$ is a quasikernel of $D$ and that $Q'\sm Q=Q'\sm V(D)$ contains every vertex $v_i$ with $v\in V(D)\sm N^+_D[Q]$ (as otherwise $Q'$ would not be a quasikernel for $D'$).  These observations imply
	\[(k+1)n-c((k+1)n)^\be \ge |Q'|\ge |Q'\sm Q|\ge k(n-|N^+_D[Q]|)\ge k(n-m),\]
	with the last step using the definition of $m$.  Rearranging this and using $n/2m\le k\le n/m$ gives
	\[c\left(\frac{n^2}{2m}\right)^\be \le c((k+1)n)^\be \le n+km\le 2n.\]
	Rearranging this yields
	\[(2m)^{\be}\ge \half c n^{2\be-1},\]
	which in turn gives \[m\ge 2^{-1-1/\be} c^{1/\be} n^{2-1/\be}\ge 2^{-3} c^{1/\be} n^{2-1/\be},\]
	giving the result.
\end{proof}

We next prove that large quasi-kernels imply small quasi-kernels.  Our proof is based  in spirit off the argument of Ai, Liu, and Peng~\cite[Proposition 4]{ai2024variable} showing large quasi-kernels (at the linear scale) imply small quasi-kernels.  Roughly, their argument works by taking a smallest quasi-kernel $Q$, decomposing it into sets $Q_1\cup Q_2$ in some way, constructing a weighted subgraph $\tilde{D}\sub D$ with weights based on how many vertices of $Q_2$ are covered, finding a large quasi-kernel $\tilde{Q}$ of $\tilde{D}$, and finally showing that if $|Q|$ is large then we can use $\tilde{Q}$ to replace $Q_2$ with a smaller set while maintaining the quasi-kernel condition.  This gives a contradiction, showing that we can not have $|Q|$ large.

Roughly, our refinement of this approach gets rid of the weighting on $\tilde{D}$ and the singular $\tilde{Q}$ and replaces it with a small set of quasi-kernels whose neighborhoods cover $\tilde{D}$ in the following sense.

\begin{lem}\label{partition by neighborhoods}
	If $0<c,\al<1$ are real numbers such that every $n$-vertex digraph $D'$ has a quasi-kernel $Q'$ with $|N^+_{D'}[Q']|\ge c n^\al$; then for every digraph $D$ there exist quasi-kernels $R_1,\ldots,R_t$ of $D$ satisfying $\bigcup_{i=1}^t N^+_{D}[R_i]=V(D)$ and $t\le \frac{1}{c(1-\al)} |V(D)|^{1-\al}$.
\end{lem}
\begin{proof}
	Let $D_1=D$.  Iteratively, if $D_{i}$ is not the empty digraph then let $Q_i\sub V(D_{i})$ be a quasi-kernel of $D_i$ satisfying $|N_{D_i}^+[Q_i]|\ge c |V(D_{i})|^\al$ and let $D_{i+1}=D_{i}-N_{D_i}^+[Q_i]$.
	
	We claim that this process stops at $D_{t+1}=\emptyset$ for some $t\le  \frac{1}{c(1-\al)} |V(D)|^{1-\al}$.  Indeed, let $n_i:=|V(D_i)|$.  By construction $n_{i+1}\le n_{i}-c n_{i}^\al$, which together with concavity implies
	\[n_{i}^{1-\al}-n_{i+1}^{1-\al}\ge (1-\al)n_{i}^{-\al}(n_{i}-n_{i+1}) \ge c(1-\al).\]
	Since the sequence $n_i^{1-\al}$ decreases by at least $c(1-\al)$ at each step and since its initial term is $|V(D)|^{1-\al}$, it follows that this process can continue for at most $\frac{1}{c(1-\al)}|V(D)|^{1-\al}$ steps as desired.
	
	Since each $Q_i$ is in particular an independent set of $D_i$, it follows that the $D_i$ digraphs satisfy the conditions of \Cref{CL-sequence}, and as such \Cref{CL-sequence}(a) implies there exist quasi-kernels $R_1,\ldots,R_t$ of $D$ with $Q_i\sub R_i$ for all $i$.  By construction
	\[\bigcup_{i=1}^t N_D^+[R_i]\supseteq \bigcup_{i=1}^t N_{D_i}^+[Q_i]=V(D),\]
	so this establishes the result.
\end{proof}

We now use this to prove large quasi-kernels imply small quasi-kernels.

\begin{lem}\label{large implies small}
	If $0<c,\al<1$ are real numbers such that every $n$-vertex digraph $D'$ has a quasi-kernel $Q'$ with $|N^+_{D'}[Q']|\ge c n^\al$; then every $n$-vertex source-free digraph $D$ has a quasi-kernel $Q$ with $|Q|\le n-\frac{c(1-\al)}{4} n^{\frac{1}{2-\al}}$.
\end{lem}
\begin{proof}
	Let $D$ be an $n$-vertex source-free digraph and let $Q\sub V(D)$ be a quasi-kernel of minimum size.
	
	\begin{claim}\label{Q1 exists}
		There exists $Q_1\sub Q$ with $N^+(Q_1)=N^+(Q)$ and with $|Q_1|\le |N^+(Q)|$.
	\end{claim}
	\begin{proof}
		By definition for each $v\in N^+(Q)$ there exists some $v'\in Q$ with $v\in N^+(v')$, and taking $Q_1=\{v':v\in N^+(Q)\}$ gives the result.
	\end{proof}
	Fix such a $Q_1$ and define the sets
	\[Q_2=Q\sm Q_1, \qquad F=V(D)\sm N^+[Q],\]
	where $F$ is meant to stand for the vertices which are ``far'' from $Q$.  \Cref{Q1 exists} gives some control over $|Q_1|$ and it remains to bound $Q_2$.  To do this, for each $R\sub F$ let
	\[C_R=\{v\in Q_2: \dist(R,v)\le 2\}.\]
	\begin{claim}
		If $R\sub F$ is a quasi-kernel of $D[F]$, then $|C_R|\le |R|$.
	\end{claim}
	\begin{proof}
		Define the set
		\[Q'=R\cup (Q_2\sm C_R) \cup (Q_1\sm N^+_D(R)).\]
		We claim $Q'$ is a quasi-kernel of $D$.  Indeed, $R,Q$ being quasi-kernels implies there are no arcs within $R$ or $(Q_2\sm C_R) \cup (Q_1\sm N^+_D(R))\sub Q$, there are no arcs from $Q_1\cup Q_2$ to $R$ by definition of $R\sub F$, and there are no arcs from $R$ to either $Q_1\sm N^+_D(R)$ or $Q_2\sm C_R$ by definition.  We have that $R$ reaches all of $F$ within two steps by definition of $R$ being a quasi-kernel of $D[F]$, that $R\cup (Q_2\sm C_R)$ reaches all of $Q_2$ within two steps by definition of $C_R$, and that $R\cup (Q_1\sm N^+_D(R))$ reaches all of $Q_1$ within one step and hence all of $N^+_D(Q)$ by definition of $Q_1$.  In total $Q'$ reaches all of $Q\cup N^+_D(Q)\cup F$ which equals $V(D)$ since $Q$ is a quasi-kernel.
		
		Now because we chose $Q$ to be a minimum quasi-kernel this implies
		\[|Q|\le |Q'|=|R|+|Q|-|C_R|-|Q_1\cap N^+_D(R)|\le |R|+|Q|-|C_R|,\]
		and rearranging gives the desired bound.  As an aside, the role of $Q'$ in this argument is analogous to the strictly smaller quasi-kernel mentioned in our comparison with the work of Ai, Liu, and Peng~\cite{ai2024variable}.
	\end{proof}
	Crucially, we observe that for every $v\in Q_2$ there exists some $u\in F$ with $v\in N^+(u)$.  Indeed because $D$ is source-free there exists some $u$ with $v\in N^+(u)$.  We can not have $u\in Q$ since this is an independent set, and if we had $u\in N^+(Q)$ then $Q\sm \{v\}$ would be a strictly smaller quasi-kernel due to \Cref{Q1 exists} guaranteeing every vertex which could be reached by $Q$ continues to be reached by $Q\sm \{v\}$, and this would contradict us defining $Q$ to be a smallest quasi-kernel.
	
	With this all in mind, let $R_1,\ldots,R_t$ be the quasi-kernels of $D[F]$ guaranteed by \Cref{partition by neighborhoods}.  Because the $N_{D[F]}^+[R_i]$ sets cover all of $F$, the observation above implies for each $v\in Q_2$ there exists some $i$ with $N_{D[F]}^+[R_i]$ containing an in-neighbor of $v$, which in turn implies $v\in C_{R_i}$.  It follows that
	\[|Q_2|\le \sum_{i=1}^t |C_{R_i}|\le \sum_{i=1}^t |R_i|\le t\cdot |F|\le \frac{1}{c(1-\al)} |F|^{2-\al}.\]
	Putting everything together, we find
	\begin{align*}|V(D)|&=|Q_1|+|Q_2|+|N^+_D(Q)|+|F|\\ &\le |N^+_D(Q)|+\frac{1}{c(1-\al)} |F|^{2-\al}+|N^+_D(Q)|+|F|\le \frac{4}{c(1-\al)}|N^+_D(Q)\cup F|^{2-\al}.\end{align*}
	Rearranging this and using $|N^+_D(Q)\cup F|=|V(D)|-|Q|$ gives
	\[|Q|\le |V(D)|-\left(\frac{c(1-\al)}{4}\right)^{\frac{1}{2-\al}}|V(D)|^{\frac{1}{2-\al}}\le |V(D)|-\frac{c(1-\al)}{4}|V(D)|^{\frac{1}{2-\al}} ,\]
	proving the result.
\end{proof}

We now put these pieces together to formally prove our equivalence result.
\begin{proof}[Proof of \Cref{equivalence}]
	Let $0<\al<1$.  First assume (a) holds, meaning there exists some $c>0$ such that every $n$-vertex digraph has a quasi-kernel $Q$ with $|N^+[Q]|\ge c n^\al$.  This assumption continues to hold if we replace $c$ with a smaller positive number, so we may assume without loss of generality that $c<1$.  It follows from \Cref{large implies small} that (b) holds with $c'=\frac{c(1-\al)}{4}>0$.
	
	Now assume (b) holds, meaning there exists some $c'>0$ such that every $n$-vertex source-free digraph has a quasi-kernel with $|Q|\le n-c' n^{\frac{1}{2-\al}}$.  Again we may assume $c'<1$, so it follows from applying \Cref{small implies large} with $\be=\frac{1}{2-\al}$ that (a) holds with $c=\frac{1}{8} (c')^{2-\al}>0$, proving the result.
\end{proof}

\section{Proof of \Cref{refinemnet}}
At a high-level, our proof showing every digraph $D$ has a quasi-kernel with $|N^+[Q]|$ large works by starting with $D_1=D$ and then iteratively defining $Q_i$ to be a large quasi-kernel of $D_i$ and taking $D_{i+1}=D_i-N^+[Q_i]$.  If we can show this process terminates after a small number of steps, then necessarily some $|N^+[Q_i]|$ must be large, which by \Cref{CL-sequence} implies there exists some quasi-kernel $Q\supseteq Q_i$ of $D$ which has $|N^+[Q]|$ large.  The simplest way to show this process takes a small number of steps would be to show each $N^+[Q_i]$ is large, but this is seemingly as hard as proving the theorem in the first place.  To get around this, we instead show $N^+[Q_i]\cap L$ is large in terms of $|L|$ for some ``nice'' set of vertices $L$, from which we will be able to lift the bound with $L$ replaced by $V(D)$ throughout.

Following Rutschmann~\cite{rutschmann2026small}, we say a vertex $v$ in a digraph $D$ is \textit{eligible} if there exists some quasi-kernel $Q$ of $D$ which contains $v$.  Let $L(D)$ denote the set of eligible vertices of $D$ and let
\[r(D)=\max_Q |N^+_D[Q]\cap L(D)|\]
where the maximum ranges over all quasi-kernels $Q$ of $D$.  Intuitively, $r(D)$ measures how many eligible vertices we can remove from $D$ by deleting some $N^+[Q]$ set like how we discussed in the sketch at the start of this section. Our main technical result gives an effective lower bound on this amount of removal.

\begin{lem}\label{eligable technical lemma}
	Every digraph $D$ has
	\[r(D)\ge \frac{1}{8}|L(D)|^{2/3}.\]
\end{lem}
\begin{proof}
	The result is trivial when $D$ is the empty digraph, so from now on we assume $D$ has at least one vertex which implies $r(D)>0$.
	
	The key idea of the proof will be to bound the number of arcs of $D$ which have both endpoints in $L(D)$, and we denote this quantity by $a_L(D)$.  We first observe
	\[a_L(D)\le r(D) |L(D)|,\]
	since by definition $|N^+[v]\cap L(D)|\le r(D)$ for any $v\in L(D)$ due to the quasi-kernel $Q$ containing $v$.  The remainder of the proof is dedicated towards lower bounding $a_L(D)$, the rough idea for which is to find a large number of disjoint sets $R\sub L(D)$ which have many arcs from vertices of $L(D)$.  For this the following will be helpful.
	
	\begin{claim}\label{lifting eligible}
		Let $D_1,D_2,\ldots$ be any sequence of digraphs with $D_1=D$ and $D_{i+1}=D_i-N^+_{D_i}[I_i]$ for some independent set $I_i$ of $D_i$.  If $R\sub V(D_i)$ is a quasi-kernel of $D_i$, then
		\[R\sub L(D),\qquad |N^+_{D}[R]\cap L(D)|\le r(D).\]
	\end{claim}
	\begin{proof}
		By \Cref{CL-sequence}(a) there exists a quasi-kernel $Q$ of $D_1=D$ which contains $R$.  This immediately implies $R\sub L(D)$ as well as
		\[|N^+_{D}[R]\cap L(D)|\le |N^+_D[Q]\cap L(D)|\le r(D).\]
	\end{proof}
	
	Note that if $|L(D)|< 4r(D)$ then the bound of the lemma trivially holds, so we may assume $|L(D)|\ge 4r(D)$. Define
	\[k=\left\lfloor \frac{|L(D)|}{4r(D)}\right \rfloor \ge 1.\]
	Let $D_1=D$.  Iteratively assume we have defined a digraph $D_i$.  Let $D_{i,1}=D_i$.  Iteratively given that we have defined $D_{i,j}$, let $R_{i,j}$ be an inclusion-maximal quasi-kernel of $D_{i,j}$ and set $D_{i,j+1}=D_{i,j}-N^+[R_{i,j}]$.  By \Cref{CL-sequence}(b) there exists a quasi-kernel $Q_i$ of $D_i$ which satisfies $R_{i,1}\cup \cdots \cup R_{i,k}\sub N^+[Q_i]$.  We then define $D_{i+1}=D_i-N^+_{D_i}[Q_i]$ and we stop this process upon defining $D_{k+1}$.
	
	Observe that each $D_{i,j}$ is a digraph of the form from \Cref{lifting eligible}.  It follows then that for each quasi-kernel $R_{i,j}$ of $D_{i,j}$ we have $R_{i,j}\sub L(D)$ and
	$|N^+_D[R_{i,j}]\cap L(D)|\le r(D)$.
	Similarly we have $|N^+_D[Q_i]\cap L(D)|\le r(D)$ for all $i$.  Thus as we construct our $D_i$ and $D_{i,j}$ digraphs by deleting either $N^+[Q_i]$ or $N^+[R_{i,j}]$, we iteratively are removing at most $r(D)$ vertices from $L(D)$, which means for all $i,j\le k$ we have
	\begin{align*}|L(D)\cap V(D_{i,j})|&\ge |L(D)\cap V(D_i)|-\sum_{s=1}^{j-1} |L(D)\cap N^+[R_{i,s}]|\\ &\ge |L(D)|-(i-1)r(D)-(j-1)r(D)\\ &\ge |L(D)|-2k r(D)+r(D)\ge  \half |L(D)|+r(D),\end{align*}
	with the second inequality using that $D_{i}$ is obtained by deleting $i-1$ sets of the form $N^+[Q_i]$ and the last inequality used $k\le \frac{|L(D)|}{4r(D)}$ by definition.
	
	Observe that each $R_{i,j}$ is a maximal independent set in $D_{i,j}$ since it is an inclusion-maximal quasi-kernel.  This in particular means every vertex in $D_{i,j}-N^+_{D}[R_{i,j}]$ has an arc to some vertex in $R_{i,j}$, and the number of such vertices lying in $L(D)$ is at least
	\[|L(D)\cap V(D_{i,j})|-|L(D)\cap N^+_{D}[R_{i,j}]|\ge (\half|L(D)|+r(D))-r(D)= \half |L(D)|.\]
	
	Thus every $R_{i,j}$ set has at least $\half|L(D)|$ arcs from vertices of $L(D)$.  Finally, we note that each of the $k^2$ sets $R_{i,j}$ are pairwise disjoint; this is evident for $R_{i,j},R_{i,j'}$, and for $R_{i,j},R_{i',j'}$ with say $i<i'$ this follows from $D_{i+1}$ being formed by deleting $N^+[Q_i]$ which contains $R_{i,j}$.  Since each $R_{i,j}$ set lies in $L(D)$, putting all this together with the previous observation implies
	\[\half k^2 |L(D)|\le a_L(D)\le r(D) |L(D)|.\]
	The definition of $k$ implies $k\ge \frac{|L(D)|}{8 r(D)}$, and plugging this in and rearranging gives \[|L(D)|^2\le 2^7 r(D)^3\le 2^9 r(D)^3.\] Taking cube roots on both sides and dividing by 8 gives the result.
\end{proof}

We conclude with our main proof of the section.

\begin{proof}[Proof of \Cref{refinemnet}]
	The upper bound on $|Q|$ for source-free digraphs will follow from applying \Cref{large implies small} with $c=\frac{1}{64}$ and $\al=\frac{2}{3}$, so it remains to prove the lower bound on $|N^+[Q]|$.
	
	Let $D_1=D$.  Iteratively given $D_i$, let $Q_i$ be a quasi-kernel of $D_i$ satisfying $|N_{D_i}^+[Q_i]\cap L(D_i)|=r(D_i)$ and define $D_{i+1}=D_i-N^+_{D_i}[Q_i]$.  Let $t$ be the smallest integer such that $D_{t+1}$ is empty.
	
	For ease of notation let $\ell_i=|L(D_i)|$ and $r_i=r(D_i)$.  By how we defined $D_{i+1}$, we have by \Cref{CL-move} that every quasi-kernel $Q'$ of $D_{i+1}$ is contained in a quasi-kernel of $D_i$, meaning $L(D_{i+1})\sub L(D_i)$.  In particular this implies $\ell_{i+1}\le \ell_i-r_i$.  Now \Cref{eligable technical lemma} implies $r_i\ge 2^{-3} \ell_i^{2/3}$ for all $i$, and using this together with concavity implies
	\[\ell_i^{1/3}-\ell_{i+1}^{1/3}\ge \frac{\ell_i-\ell_{i+1}}{3 \ell_i^{2/3}}\ge 2^{-5}.\]
	Since the sequence $\ell_i^{1/3}$ decreases by at least $2^{-5}$ and since $\ell_i\ge 0$ for all $i$, we conclude that this process must stop at 
	\[t\le 2^5 \ell_1^{1/3}.\]
	
	Now let $m$ be the maximum size of $|N^+[Q]|$ amongst all quasi-kernels $Q$ of $D$.  Because the $D_i$ are of the form from \Cref{CL-sequence}, for every quasi-kernel $Q_i$ of $D_i$ there exists a quasi-kernel $Q_i'\supseteq Q_i$ of $D_1=D$ which necessarily satisfies
	\[|N^+_{D_i}[Q_i]|\le |N^+_D[Q'_i]|\le m.\]
	Since the $N^+_{D_i}[Q_i]$ sets partition $V(D)$ by construction, we necessarily have $|V(D)|\le mt$.  We also observe by definition of $r$ and $m$ that $r(D)\le m$, so \Cref{eligable technical lemma} implies
	\[|L(D)|\le 8^{3/2} m^{3/2}.\]
	Putting all this together gives
	\[|V(D)|\le m\cdot t\le m\cdot 2^5 |L(D)|^{1/3}\le m\cdot 2^7 m^{1/2}\le  2^9 m^{3/2},\]
	and rearranging gives $m\ge 2^{-6} |V(D)|^{2/3}$, proving the result.   
\end{proof}

\bibliographystyle{abbrv}
\bibliography{QuasiRefs}
\end{document}